\batchmode
\documentclass[11pt]{amsart}
 \usepackage{amssymb}
 \usepackage{amsmath}
\usepackage{amscd} 
\begin{document}
    \newtheorem{Theorem}{Theorem}[section]
    \newtheorem{Proposition}[Theorem]{Proposition}
    \newtheorem{Lemma}[Theorem]{Lemma}
    \newtheorem{Corollary}[Theorem]{Corollary}
   \newcommand{\ra}{\rightarrow}   
   \newcommand{\ul}{\underline}
   \newcommand{\varp}{\varphi} 
   \newcommand{\al}{\alpha}
   \newcommand{\bt}{\beta}
   \newcommand{\ve}{\varepsilon}
    \newcommand{\aut}{{\rm Aut}} 
     \newcommand{\Z}{\boldsymbol{Z}} 
     \newcommand{\C}{\boldsymbol{C}} 
     \newcommand{\R}{\boldsymbol{R}} 
     \newcommand{\B}{\boldsymbol{P}} 
     \newcommand{\Q}{\boldsymbol{Q}}

\title[Non-upper-semiconitnuity of algebraic dimension]{Non-upper-semicontinuity of algebraic dimension for families of compact complex manifolds } 

\author{A. Fujiki and M. Pontecorvo} 
\footnote{Research partially supported by Grant-in-Aid for Scientific Research (No.18340017) MEXT, 
and PRIN, Italy: Metriche riemanniane e variet\'a differenziabili} 

\maketitle

\vspace{5 mm} 
\begin{abstract} In this note we show that in a certain subfamily of 
the Kuranishi family of any half Inoue surface the algebraic dimensions 
of the fibers jump downwards at special points of the parameter space 
showing that the upper semi-continuity of algebraic dimensions in any sense 
does not hold in general for families of compact non-K\"ahler manifolds.  
In the K\"ahler case, the upper semi-continuity always holds true in a certain weak sense. 
\end{abstract}

\section{Statement of results} 

In this note we shall show the following theorem 
which gives examples of holomorphic families of compact complex surfaces 
whose algebraic dimensions jump downwards under specializations.  

\begin{Theorem}\label{main} 
Let  $S$  be a half Inoue surface 
with second Betti number $m$ and  $C$  the unique twisted anti-canonical curve on it.  
Let  $g: ({\mathcal S},{\mathcal C}) \ra  T,\ (S_o,C_o)=(S,C),\ o\in T$, 
be the Kuranishi family of deformations of the pair  $(S,C)$.  
Then the Kuranishi space  $T$  is smooth of dimension  $m$  and 
contains a divisor with normal crossings  $A=\bigcup _{i=1}^m A_i$ with 
$m$ smooth irreducible components passing through the base point  $o$ such that 
the following hold: the fiber  $S_t, t\neq o$, is a blown-up half Inoue surface if $t \in  A$, 
and is a blown-up elliptic diagonal Hopf surface if $t \notin A$.  
\end{Theorem} 

We refer to Proposition \ref{mlt} below for more precise structure 
of the surface  $S_t$ for  $t\notin A$. 
Recall that the {\em algebraic dimension} of a compact connected complex surface 
is the transcendence degree of its meromorphic function field.  
It takes one of the values $0, 1$ and $2$. 
Any half Inoue surface has algebraic dimension zero, while 
any elliptic diagonal Hopf surface is of algebraic dimension one. 
This is the property of our main interest in this note: 

\begin{Corollary}\label{mainc}
Let $a_t:=a(S_t)$ be the algebraic dimension  of $S_t$.  
Then  $a_t  = 0$ for $t\in A$ and  $= 1$ for $t\not\in A$.  
Thus the algebraic dimension jumps downwards at special points; in particular 
it is not upper semicontinuous in the parameter  $t$. 
\end{Corollary}

\vspace{3 mm} 
After some preliminaries in Section 2, 
we prove Theorem \ref{main} in Section 3 together with the following supplement 
to it, which reveals special features of our examples. 

\begin{Proposition}\label{mlt} 
Let the notations be as in Theorem \ref{main}. 
Then for any $t\notin A$, $S_t$ is the blowing-up 
of an elliptic diagonal Hopf surface $\bar{S}_t$. 
The elliptic fibration of $\bar{S}_t$ is smooth 
except for two multiple fibers with multiplicity two, 
and the image $\bar{C}_t$ of $C_t$  in  $\bar{S}_t$ 
is one of its smooth fibers. 
Moreover, $S_t$  is obtained from $\bar{S}_t$ by 
blowning up $m$ points on  $\bar{C}_t$.   
\end{Proposition} 

\noindent{\bf Example}.  When $m=1$, the half Inoue surface $S$ is unique up to isomorphisms.  
It contains a unique curve $C$,  
which is a rational curve with a single node.  
We get the Kuranishi family of deformations $\{(S_t,C_t)\}_{t\in D} $ of $(S,C)$ 
parametrized by a one dimensional disc  $D=\{|t|<\ve \}, \ve>0$.  
For any $t\neq 0$,  $S_t$ is the blowing-up at one point $p_t$ 
of an elliptic diagonal Hopf surface $\bar{S}_t$ 
over the complex projective line $\B$.  
Let $\bar{C}_t$ be the fiber of  $\bar{h}_t :\bar{S}_t\ra \B$ containing $p_t$.  
Then $\bar{h}_t$ is smooth except for two multiple fibers with multiplicity two and 
$C_t$ is the proper transform of $\bar{C}_t$ in  $S_t$.  
In particular  $a_t=1$ for $t\neq 0$ and $a_0 = 0$ (cf.\ Proposition \ref{mlt} below).  

\vspace{3 mm} 
{\em Remark 1.1}. For each integer  $n>1$  
by taking the product of $n$ copies $S_t^n$ of each fiber of  $f$  in Theorem 1.1 or 
by considering the Douady space (Hilbert scheme) $S_t^{[n]}$ of $0$-dimensional 
analytic subspaces of  $S_t$  of length  $n$  we get a family of $2n$-dimensional 
compact complex manifolds $X_t$ with the same parameter space $T$ 
such that $a(X_t)=0$ for $t\in A$, but $a(X_t)=n$ for $t\notin A$.

\vspace{3 mm} 
In the above results the surfaces $S_t$ are all non-K\"ahler surfaces. 
In fact the phenomena 
as above never occur in a family of compact K\"ahler manifolds. 
In Section 4, for the purpose of comparison, 
we give a general property of the variation of algebraic dimensions of the fibers 
in a family of compact K\"ahler manifolds. 

\section{Preliminaries}  

A {\em cycle of rational curves} on a smooth surface is a compact connected curve 
$C$  which is either an irreducible rational curve with a single node or 
is a reducible curve with $k$ nodes 
whose irreducible components are nonsingular rational curves 
$C_i, 1\leq i\leq k, k\geq 2$, such that $C_i$ and $C_{i+1}$ intersects at a single point 
and there exists no other intersections, where $C_{k+1}=C_1$ by convention. 

A compact connected complex surface $S$ is called of {\em class VII} if 
$b_1(S)=1$ and $\kappa (S)=-\infty$, 
where $b_1$ and $\kappa $ are the first Betti number and 
the Kodaira dimension of $S$ respectively.  
Such a surface is non-K\"ahler and its algebraic dimension is less than two. 
$S$ is called of {\em class VII$^+_0$} 
if $S$ is of class VII, is minimal, i.e., contains no $(-1)$-curves, and 
with positive second Betti number. 
(A $(-1)$-curve is a nonsingular rational curve with self-intersection number $-1$.)

Inoue \cite{ino1}\cite{ino} constructed the first examples of surfaces 
of class VII$^+_0$.  
These surfaces are called hyperbolic, parabolic and half Inoue surfaces 
according to Nakamura \cite{na84}.  
By \cite[(8.1)]{na84} 
a surface $S$  of class VII$^+_0$ is a hyperbolic (resp.\ parabolic) Inoue surface
if and only if   $S$ 
contains two cycles of rational curves $C_1$ and $C_2$ 
(resp.\ a cycle of rational curves  $C_0$ and a smooth elliptic curve $E$). 
Similarly, 
$S$ is a half Inoue surface if and only if 
$S$ contains a cycle $C$ of rational curves 
and there exists an unramified double covering $u: \tilde{S} \ra  S$ such that 
$u^{-1}(C)$ is a disjoint union of two cycles 
of rational curves $\tilde{C}_1$ and $\tilde{C}_2$ 
which are mapped isomorphically onto  $C$ (cf.\ \cite[(1.6)(9.2.4)]{na84}).  
$\tilde{S}$ is then a hyperbolic Inoue surface.  
All these Inoue surfaces have infinite cyclic fundamental group.  
By \cite{ino} the algebraic dimension of these surfaces equals zero. 

In the hyperbolic (resp.\ parabolic) case 
denote by  $C$  the union of $C_1$ and $C_2$ (resp.\ $C_0$ and $E$). 
Then for any Inoue surface 
there exists no irreducible curves other than the irreducible components 
of  $C$.  Moreover, in the hyperbolic or parabolic Inoue case  
$C$  is the unique {\em anti-canonical curve} on  $S$, i.e., 
$C$ is the unique member of the anti-canonical system $|-K_S|$. 
Similalry, in the half Inoue case  $C$  is the unique $L$-{\em twisted anti-canonical 
curve} in the sense that  $C$  is the unique member of the $L$-twisted anti-canonical 
system $|-(K_S+L)|$ for a unique non-trivial holomorphic line bundle $L$ with 
$L^{\otimes 2}$ trivial.  

A {\em diagonal Hopf surface} is a Hopf surface $S=S(\al,\bt)$ obtained as the quotient 
of $\C^2-0$ by an infinite cyclic group generated by the diagonal transformation 
$(z,w) \ra (\al z,\bt w)$ for complex numbers $\al , \bt$ with  $0<|\al |,|\bt |<1$. 
The images of $z$- and $w$- axes give two canonical elliptic curves $E_1$ and $E_2$ 
on  $S$  and we have  $E_1+E_2=-K_S$ (cf.\ \cite[(97)]{kd2}).  
If $S$ admits a non-constant meromorphic function, 
it admits a unique structure 
of an elliptic surface  $h: S \ra \B$  over the projective line  $\B$ and 
any irreducible curve on  $S$  is smooth and is a fiber of  $h$  
up to multiplicity (cf.\ \cite{kd2}).  
Recall also that $S$ has the vanishing second Betti number. 

We also quote two lemmas \cite[Lemmas 3.2, 3.3]{fp} for later purpose 
for the convenience of the reader.   
The first one is originally due to Nakamura and the second one easily 
follows from the first one.  

\begin{Lemma}\label{aatc}
Let  $\hat{S}$  be a compact complex surface of class VII$_0$ 
with infinite cyclic fundamental group. 
Suppose that $\hat{S}$ admits a disconnected anti-canonical curve. 
Then  $\hat{S}$  is either a hyperbolic or parabolic Inoue surface 
or a diagonal Hopf surface. 
\end{Lemma} 

\begin{Lemma}\label{atc} 
Let  $\tilde{S}$  be a compact complex surface of class VII 
with infinite cyclic fundamental group. 
Suppose that $\tilde{S}$ admits a disconnected anti-canonical curve  $\tilde{C}$.  
Then the minimal model  $\hat{S}$  of  $\tilde{S}$  is 
either a hyperbolic or parabolic Inoue surface 
or a diagonal Hopf surface, and  $\tilde{S} \ra \hat{S}$  is obtained by blowing up $\hat{S}$ 
at a finite number of points (possibly infinitely near) 
on the image $\hat{C}$ of $\tilde{C}$. 
Moreover, $\hat{C}$ is an anti-canonical curve  on  $\hat{S}$. 
\end{Lemma} 

\vspace{3 mm} 
{\em Remark 2.1}. 
It is known that the fundamental group
of a known surface $S\in {\rm VII}_0$ is always infinite cyclic 
unless $S$ is a Hopf surface; in the latter case however 
$S$ is always finetely covered by a primary Hopf surface,  
that is one with $\pi_1(S)=\Z$. 

\vspace{3 mm} 

\section{Proof of main results}  

We consider the Kuranishi family 
$g: ({\mathcal S},{\mathcal C}) \ra  T,\ (S_o,C_o)=(S,C),\ o\in T$ 
of deformations of the pair $(S,C)$ as in Theorem \ref{main}. 
This was studied 
in our previous paper \cite{fp}, to which we refer for the details. 
Indeed, 
except the assertion that $S_t, t\notin A$, has the structure of an 
elliptic surfaces with special structures, all statements in the theorem 
are given in Proposition 3.14 of \cite{fp} partly without proof,  
as an analogue of the corresponding result 
in the hyperbolic case (\cite[Proposition 3.13]{fp}). 

We first prove Lemma \ref{atch} below which is a detailed version of 
Lemma 3.4 of \cite{fp} also stated without proof there.  
Let  $S$  in general be a compact complex surface of class VII 
with infinite cyclic fundamental group.  
Then there exist a unique unramified double covering $u: \tilde{S} \ra S$ and 
a unique non-trivial holomorphic line bundle $L=L_S$ 
with $L^{\otimes 2}$ trivial such that $u^*L$ is trivial.  
In fact they are both associated 
to the representation of the fundamental group 
$\pi_1(S)\cong \Z \ra \Z_2=\{\pm 1\}\subseteq \C^*$. 

\begin{Lemma}\label{atch}  Let the notatios be as above. 
Suppose that $S$ contains an $L$-twisted connected anti-canonical curve $C$ 
such that $\tilde{C}:=u^{-1}(C)$ is disconnected in  $\tilde{S}$. 
Then the minimal model  $\bar{S}$  of  $S$  is either a half Inoue surface 
or an elliptic diagonal Hopf surface.  
In the latter case 
the elliptic fibration of $\bar{S}$ 
is smooth except for two multiple fibers with multiplicity two, 
and the image $\bar{C}$ of  $C$  in  $\bar{S}$ 
is one of the smooth fibers. 
$S$  is obtained by blowning up points on $\bar{C}$ (unless $S=\bar{S}$) and 
$C$  is the proper transform of $\bar{C}$.  
\end{Lemma}

{\em Proof}.  Let  $\hat{S}$ be the minimal model of $\tilde{S}$. 
Let $\tilde{v}: \tilde{S} \ra  \hat{S}$ and 
$v: S \ra  \bar{S}$ be the blowing-down maps. 
For any $(-1)$-curve $B$ on  $S$,  
$u^{-1}(B)$ is a disjoint union of two $(-1)$-curves 
$\tilde{B}_1$ and $\tilde{B}_2$ on  $\tilde{S}$.  
Thus contracting successively all the $(-1)$-curves on  $\tilde{S}$  obtained 
in this way we get a blowing-down map $\tilde{v}': \tilde{S} \ra  \hat{S}'$.  
We show that $\hat{S}'$ coincides with  $\hat{S}$.  
Let  $\iota $ be the Galois involution for  $u$. 
Suppose that there exists a $(-1)$-curve $\tilde{B}$ on  $\hat{S}'$. 
$\iota (\tilde{B})$ is a $(-1)$-curve, 
but by the definition of $\hat{S}'$ it must intersect with 
$\tilde{B}$ so that $D:=\tilde{B}\cup \iota (\tilde{B})$ is connected 
with intersection number $\tilde{B}\cdot\iota(\tilde{B})\geq 2$. 
Then $D^2\geq 2$, which is impossible since $a(\hat{S}')=a(S)\leq 1$ 
(cf.\ \cite[Th.8]{kd5}). 
Hence  $\hat{S}'=\hat{S}$ and we get 
the induced double covering  $\hat{u}: \hat{S} \ra  \bar{S}$ with Galois involution 
denoted by $\hat{\iota}$.  

On the other hand, since $C=-(K_{S}+L)$ by our assumption, 
we have $\tilde{C}=-K_{\tilde{S}}$. 
$\tilde{C}$ is thus an anti-canonical curve on  $\tilde{S}$, and 
hence by Lemma \ref{atc} 
the image $\hat{C}$ of $\tilde{C}$ in  $\hat{S}$ 
under $\tilde{v}: \tilde{S}\ra \hat{S}$ is 
again a disconnected anti-canonical curve on  $\hat{S}$ 
and $\tilde{v}$ is obtained by blowingn-up points on $\hat{C}$.  
Moreover, each of the two connected components $\hat{C}_\al, \al=1,2$,  
of  $\hat{C}$ is mapped isomorphically onto the image  $\bar{C}$ on $\bar{S}$. 
Thus by Lemma \ref{aatc} $\hat{S}$ is 
either a hyperbolic Inoue surface or a diagonal Hopf surface 
and we have $-K_{\hat{S}}=\hat{C}_1+\hat{C}_2$.  
(Since $\hat{C}_1 \cong \hat{C}_2$,  $\hat{S}$ is never 
a parabolic Inoue surface.)  
Therefore in the first case  $\bar{S}$  is a half Inoue surface, 
and in the second case it is a Hopf surface with infinite cyclic fundamental group, 
i.e., a primary Hopf surface. 

We now consider the second case in more detail.  
There are two types of primary Hopf surfaces as exhibited in \cite[(94)(95)]{kd2}, 
of which one type consists precisely of diagonal Hopf surfaces.  
These two types are preserved under a finite unramified covering, as 
follows easily from their defining formulae (loc.cit).  
Thus  $\bar{S}$ is a diagonal Hopf surface as well as $\hat{S}$,  
and the first assertion is proved.  Note also that if $\bar{S}=S(\al,\bt)$, 
then $\hat{S}$ is identified with $S(\al^2,\bt^2)$.  

We further show that  $\bar{S}$ is an elliptic surface. 
Note first that each $\hat{C}_\al$ is a smooth elliptic curve as well as $\bar{C}$. 
On the other hand, 
since $\bar{S}$ is a diagonal Hopf surface, 
$\bar{S}$ admits two elliptic curves $\bar{E}_1$ and $\bar{E}_2$ 
such that $-K_{\bar{S}}= \bar{E}_1+ \bar{E}_2$. 
Moreover, if $\hat{E}_i:=\hat{u}^{-1}(\bar{E}_i), i=1,2$, 
these two curves are again the canonical elliptic curves on 
the diagonal Hopf surface $\hat{S}$, and hence we have 
$\hat{E}_1+\hat{E}_2=-K_{\hat{S}}=\hat{C}_1+\hat{C}_2$.   
Note that 
$\hat{E}_1+\hat{E}_2\neq \hat{C}_1+\hat{C}_2$ since their images 
$\bar{E}_1+\bar{E}_2$ and $\bar{C}$ in  $\bar{S}$ are different.  
Thus the anti-canonical system $|-K_{\hat{S}}|$ has positive dimension and 
so the algebraic dimension of $\hat{S}$ is positive. 
Hence $\hat{S}$ admits a unique elliptic fibration  $\hat{h}: \hat{S} \ra  \B$ 
such that the above curves $\hat{E}_i, i=1,2$, and $\hat{C}_\al ,\al=1,2$,  
become fibers of  $\hat{h}$. 

On the other hand, since the arithmetic genus  $\chi (O_{\hat{S}})$ of $\hat{S}$ 
vanishes, from the canonical bundle formula 
for elliptic surfaces \cite[Th.12]{kd5}, we deduce easily \begin{equation}\label{aa} 
-\hat{K}=\hat{h}^*O(2)-\sum_\nu (m_\nu -1)F_\nu , 
\end{equation}  
where $O(l)$ is the line bundle of degree $l$ on  $\B$, 
and $F_\nu ,0\leq \nu \leq b$, are multiple fibers with multiplicity $m_\nu >1$. 
(In the case of a diagonal Hopf surface we have $0\leq b\leq 2$ by \cite[Th.31]{kd2}.) 
Thus, if  $b=2$, noting that $m_\nu F_\nu = \hat{h}^*O(1)$ 
we have  $-\hat{K}=F_1+F_2$.  Since $F_i$ do not move in $\hat{S}$, 
this implies that dim $|-\hat{K}|=0$, contradicting what we have obtained.  
Thus $b\neq 2$.  

Next we show that $b\neq 1$. 
In fact, suppose that $b=1$ so that $-\hat{K}=\hat{h}^*O(1)+F_1$. 
Then neither of  $\hat{C}_\al$  is a multiple fiber since the multiple fiber 
is fixed by $\hat{\iota} $ while $\hat{\iota}(\hat{C}_1)=\hat{C}_2$.  
Thus $\hat{C}_\al = \hat{h}^*O(1)$  and hence 
 $-\hat{K}=\hat{C}_1+\hat{C}_2=\hat{h}^*O(2)$, 
which is a contradiction to (\ref{aa}) since $b=1$.  
Thus $b=0$ and $\hat{h}$ is a principal elliptic bundle.  

Now $\bar{\iota}$  induces a non-trivial involution on  $\B$ 
making $\hat{h}$ equivariant, and 
we have the induced morphism 
$\bar{h}: \bar{S}\cong \hat{S}/\langle \iota\rangle \ra 
\B/\langle \iota\rangle \cong \B$ 
giving the elliptic fibering structure on  $\bar{S}$.  
Thus $\bar{h}$  has exactly two multiple fibers of multiplicity two over 
the two fixed points of $\bar{\iota}$ on  $\B$ and is otherwise smooth.  

Finally, each of the inverse images in $\hat{S}$ 
of the two multiple fibers of $\bar{h}$ 
is $\hat{\iota}$-invariant, 
while we have $\hat{\iota}(\hat{C}_1)=\hat{C}_2$.  
Thus  $\bar{C}$ is a smooth fiber of $\bar{h}$.  
The last assertion immediately follows from the corresponding assertion 
for  $\tilde{v}$ and $\hat{C}$ proved above. 
\hfill $\square$  

\vspace{3 mm} 

{\em Proof of Theorem \ref{main} and Proposition \ref{mlt}}.  
That  $T$  is smooth of dimension $m$ was shown in Proposition 3.12 of 
\cite{fp}.  
The hypersufaces $A_i,1\leq i\leq m$, of the theorem are $T(p)$ in 
\cite[Prop.3.12 ]{fp} defined for each node $p$ of $C$.   
In fact $C$  is a cycle of rational curves with $m$ irreducible components. 
So it admits precisely $m$ nodes.  
$T(p)$ is the locus of the point $t\in T$  such that  
the node $p$ remains to be a node in $C_t$.  
(That is, along $T(p)$ the deformation of the 
isolated singularity germ $(C,p)$ is trivial.)  
Thus  $C_t$ is still a cycle of rational curves if $t\in A$ and 
is a nonsingular elliptic curve if $t\notin  A$.  

Now we take the unique unramified double covering  $u: \tilde{S} \ra  S$ over $S$ and 
let $\tilde{C}:=u^{-1}(C)$.  $\tilde{C}$ has two connected components 
each of which is mapped isomorphically onto  $C$. 
We can extend this covering to a relative double covering of the family so that 
for each  $t$ we have the unramified double covering  $u_t: \tilde{S}_t \ra  S_t$. 
Then similarly $\tilde{C}_t:=u_t^{-1}(C_t)$ has two connected components,  
each mapped isomorphically onto  $C_t$.  Moreover, 
by Proposition 3.13 of \cite{fp} 
$\tilde{C}_t$ is the 
anti-canonical curve on  $\tilde{S}_t$. 
Let $L_t$ be the unique non-trivial holomorphic line bundle on  $S_t$ 
with $L_t^{\otimes 2}$ trivial, i.e., 
$L_t=L_{S_t}$, which depends holomorphically on $t$ with $L_o=L$.  
 
$C_t$ is then an $L_t$-twsited anti-canonical curve on $S_t$. 
Indeed, $u_t^*(K_t+C_t)=\tilde{K}_t+\tilde{C}_t=0$ on $\tilde{S}_t$. 
Thus we get $K_t+C_t=L_t$ or $0$. 
But since $L_o=L$, by continuity the former holds true for any $t$ 
as desired. 
Then by Lemma \ref{atch} the minimal model $\bar{S}_t$ of $S_t$ is 
either a half Inoue surface or a diagonal elliptic Hopf surface. 
By the above description of the curve $C_t$ these two cases occur 
precisely when  $t\in A$ and $t\notin A$ respectively. 
Together with Lemma \ref{atch} 
the theorem and the proposition follow.  
(If $t\notin A$, $b_2(\bar{S}_t)=0$ and $b_2(S_t)=m$. Thus the blowing-up ocurrs at 
$m$ points on  $\bar{C}_t\subseteq \bar{S}_t$.) 
\hfill $\square$ 

\vspace{3 mm} 
{\em Remark 3.1}.   
Lemma 3.4 and Proposition 3.14 in \cite{fp} follow clearly from 
Proposition \ref{mlt} and Lemma \ref{atch} above.  
In fact, the arguments above are completely the same even if we start from 
a properly blown-up half Inoue surface instead of just a half Inoue surface 
as we did in \cite{fp}. 

\vspace{3 mm} 
{\em Remark 3.2}. 
1)   The above results can also be stated in terms of the Kuranishi family 
\begin{equation}\label{kur2} 
\tilde{g}: (\tilde{S},\tilde{C}) \ra  \tilde{T},\ 
(\tilde{S}_o,\tilde{C}_o)=(\tilde{S},\tilde{C}),\ o\in \tilde{T}, 
\end{equation}  
of deformations of the pair  $(\tilde{S},\tilde{C})$, 
where, as in the above proof, 
$u: \tilde{S} \ra  S$  is the canonical unramified double covering of  $S$  
with the Galois involution $\iota $ on  $\tilde{S}$ and $\tilde{C}=u^{-1}(C)$. 
$\tilde{T}$  is known to be smooth of dimension $2m$ and the family is universal 
(\cite[Prop.3.13]{fp}).  
Thus $\iota $ extends to an involution $\tilde{\iota}$ on the family  $\tilde{g}$.  
Let  $\tilde{T}^\iota $ be the fixed point set of $\tilde{T}$.   
Then the restriction of this family to  $\tilde{T}^\iota $ gives 
the universal family of deformations 
of the triple $(\tilde{S},\tilde{C},\iota )$, and 
the (fiberwise) quotient by $\tilde{\iota }$ 
of the restricted family to $\tilde{T}^\iota $ 
is naturally identified with the Kuranishi family  $g$ 
of the pair  $(S,C)$ considered above.  
Thus the same statement as Corollary \ref{mainc} holds true also for the 
family over $\tilde{T}^\iota $ above.  

It may be interesting to ask if such a phenomenon is accidental or 
a rather general one.  For instance 
we may ask the following question: 
Starting from any hyperbolic Inoue surface $\tilde{S}$ (which is in general 
not an unramified double covering of a half Inoue surface), 
consider the Kuranishi family (\ref{kur2}) of deformations of $(\tilde{S},\tilde{C})$, 
where $\tilde{C}$ is the unique anti-canonical curve on  $\tilde{S}$. 
Then can we find a subspace  $\tilde{T}'$  of  $\tilde{T}$ 
with $o \in  \tilde{T}'$ such that when the family is restricted to $\tilde{T}'$ 
we get the non-upper-semicontinuity of algebraic dimension as in Corollary \ref{mainc} ? 
How can we characterize geometrically such families ? 

2)  Some holomorphic line bundles defined on  $g^{-1}(T-A)$  do not extend 
to the total space ${\mathcal S}$ of the Kuranishi family of Theorem \ref{main}. 
Let us discuss 
this phenomenon in the simplest case $m=1$ (cf.\ Example in Section 1). 
In this case for $t\neq o$, $S_t$ has a canonical structure of an elliptic surface 
$f_t: S_t \ra \B$. We consider the line bundle  $H_t:=f_t^*O(1)$ on  $S_t$.  
$C_t$ is an irreducible component of a fiber, say $F_b$,  of  $f_t$ 
over the point  $b=b_t\in \B$; 
in fact $F_b=C_t+E_t$ for a unique $(-1)$-curve $E_t$ on  $S_t$. 
Obviouly, $K_t, L_t$ and 
$C_t,t\neq o$, extend holomorphically to  $S=S_o$, while $F_b$ and $E_t$ do not. 
More precisely we have 

\vspace{2 mm} 
{\bf Claim}.  The line bundles $F_t$ and $[E]_t, t\neq o$, never extend holomorphically 
to a line bundle on  $S=S_o$, where  $[E]_t$ is the line bundle defined by  $E_t$.  

\vspace{2 mm} 
{\em Proof}.  Suppose that  $F_t$ extends holomorphically to $F_o$ on  $S$.  
Then by the upper semicontinuity we have $2=h^0(F_t)\leq h^0(F_o)$, where 
$h^0(F_s)=\dim H^0(F_s), s\in D$.  
Since $a(S)=0$, we must have  $h^0(F_o)\leq 1$, which is a contradiction. 
The result for $[E]_t$ then follows from the relation $F_b=C_t+E_t$. Or similarly 
to the above, the extension implies that  $1=h^0([E]_t)\leq h^0([E]_o)\leq 1$.  
Then a nonzero section of $h^0([E]_o)$ would give a curve other than  $C$ on $S$, 
which is a contradiction.  \hfill $\square$ 

\vspace{3 mm} It seems intersting to study the behaviour of  $E_t$ and $[E]_t$ 
when $t$ tends to $o$.  

\section{K\"ahler case} 

The phenomenon as in Theorem \ref{main} never occurs in the K\"ahler category. 
We shall give in this section a proof of this fact 
known certainly to experts for the convenience of the reader.  Namely 
we prove: 
\begin{Proposition}\label{ald}
Let  $f: X \ra  T$  be a smooth holomorphic family of compact connected K\"ahler manifolds  
parametrized by a connected complex manifold $T$.  
Suppose that  $T$ is simply connected or dim $T=1$.  
For any nonnegative integer $k$  let  $T_k$ be the subset of  $T$  
consisting of points $t \in  T$  
such that the fiber $X_t:=f^{-1}(t)$ is of algebraic dimension $a(X_t)\geq k$.   
Then 
$T_k$ is the union of at most a countable number of analytic subsets of $T$.   
\end{Proposition} 
In general the sets $T_k, k>0$, are not closed so that algebraic dimension is 
not upper semi-continuous with respect to the classical topology even 
in the K\"ahler case.  $T_k$ can actually be dense as the following examples show. 

\vspace{3 mm} 
\noindent {\bf Example}.  
1) Let  $S$  be a K3 surface or a complex torus of dimension two.  
Let  $f: {\mathcal S} \ra T, o\in T, S_o=S$, be the Kuranishi family of  $S$, where 
$T$ is a smooth germ of dimension 20 and 4 respectively.  
Then the locus $T^+:=T_1\cup T_2$ in $T$ where the fiber is of positive algebraic dimension 
is a countable union of smooth hypersufaces.  Moreover, $T^+$ is dense in $T$.

2) Let $S$ be as above.  For each K\"ahler class on $S$ we have the associated 
twistor space  $Z$ which is a smooth fiber space over the complex projective 
line $\B^1$.  The locus $B^+$ of $\B$  where the fiber is of positive algebraic dimension is 
a countable dense subset of  $\B$.  For an explicit example of the set $T^+$ see e.g. 
\cite[Prop.5.7]{fjt}. 

\vspace{3 mm} 
In 1) and 2) above the points corresponding to projective surfaces are even dense 
in the base space.  This is a consequence  of \cite[Th.4.8]{fj83}, which in 
general holds for any compact hyperk\"ahler manifolds.  
Similar phenomena are also common in the non-K\"ahler case. For instance 
the diagonal Hopf surfaces are parametrized by $(\al ,\bt )\in D^*\times D^*=:T$ with the 
corresponding surface $S(\al ,\bt )$ given by $(\C^2-\{0\})/\langle (\al ,\bt )\rangle $, 
where $D^*$ is the punctured unit disc.  By \cite[Th.31]{kd2} $a(S(\al ,\bt ))=1$ if and only if 
$\al ^m=\bt ^n$ for some $(m,n)\in \Z^2$.  $T^+$ then forms a dense counable family of 
analytic subsets of $D^*\times D^*$.  

\vspace{3 mm} 
Now for the proof of Proposition \ref{ald} 
we first note that for any compact connected complex manifold  $X$  
the algebraic dimension $a(X)$ is described as 
\begin{equation}\label{a} 
 a(X)= \mbox{max} \{\kappa (L); L\in Pic X\}  
\end{equation}  
where  $\kappa (L)$ is the Iitaka (or $L$-) dimension of the line bundle $L$ and 
$Pic X$ is the Picard variety of  $X$ (cf.\ \cite{u}). 
Next, we recall the following result of Lieberman and Sernesi \cite{ls}. 

\begin{Lemma}\label{ls} 
Let  $f: Y \ra  S$  be a proper smooth morphism of irreducible complex spaces 
with connected fibers. 
Let  $L$  be a holomorphic line bundle on  $Y$.  
For any non-negative integer  $k$ 
we set  $S_k(L):=\{s\in S; \kappa (L_s)\geq k\}$, 
where $L_s$ is the restriction of $L$ to the fiber  $Y_s$. 
Let $k_0:= \min\{k; S_k(L)\neq S\}$.   
Then $S_{k_0}(L)$ is the union of at most a countable number 
of analytic subvarieties of  $S$.  
\end{Lemma} 

By applying the above lemma to each of the subvarieties of which $S_{k_0}(L)$ is 
a union, we get that  $S_{k_0+1}(L)$ again is 
the union of at most a countable number of proper analytic 
subvarieties of  $S$.  Proceeding in the same way, we obtain the following: 

\begin{Corollary}\label{lss} 
For each $k\geq 0$,  
$S_k(L)$  is 
the union of at most countably many analytic subvarieties of  $S$.  
\end{Corollary} 

{\em Proof of Proposition \ref{ald}}.  
Let $f: X \ra  T$ be as in the proposition.  
Suppose first that $T$  is simply connected. 
Then the local system  $R^2f_*\Z$ is trivial and is identified with the trivial system 
$T\times \Gamma  \ra  T$, where  $\Gamma :=H^2(X_o,\Z)$ for some fixed reference point $o\in T$. 
Any element $\gamma $ of $\Gamma $ thus defines 
a (constant) section  $s(\gamma )$ of $R^2f_*\Z$. 
We consider the long exact sequence 
\[	 \ra  R^1f_*O^*_X  \ra  R^2f_*\Z  \ra  R^2f_*O_X  \ra  \]
obtained by applying $Rf_*$ to the exponential sequence on $X$ 
\[   0 \ra  \Z \ra  O_{X} \ra  O^*_{X} \ra  0 . \]
For any  $\gamma \in \Gamma $, 
$s(\gamma )$ is mapped to a section $b(\gamma )$ of $R^2f_*O_X$. 
Since  $X_t$ are all K\"ahler, $R^2f_*O_X$ is a locally free $O_T$-module 
and 
hence the zero locus of  $b(\gamma )$ is an analytic subset $T_\gamma $ of $T$. 
Also by the K\"ahler assumption, $h^i(O_{X_t})$ is constant for all $i$ 
and therefore $O_X$ is cohomologically flat in all dimensions with respect to $f$ 
\cite[Th.4.1.2]{bs}. 
In particular  $T_\gamma $ coincides with the set of points $t\in T$ 
such that the restriction $s(\gamma )_t\in H^2(X_t,\Z)$ is of type 
$(1,1)$ on  $X_t$, or equivalently, 
$s(\gamma )_t$ is the first chern class of a holomorphic line bundle on  $X_t$.  

Let  $p: P=Pic X/T  \ra  T$  be the relative Picard variety associated 
to the morphism  $f$ (cf.\ the remarks at the end of \cite[\S 3]{gr} and \cite{bi}).  
Let  $p_\gamma: P_\gamma  \ra  T$ be the component of $P$  
parametrizing line bundles with chern class $\gamma$.   
Then the image of  $p_\gamma $ is precisely the set $T_\gamma $, and 
$p_\gamma $ is proper, smooth and with connected fibers over  $T_\gamma $ again 
by the K\"ahler condition. 
(Here we consider $T_\gamma $ as a reduced complex space.) 
Take the fibered product  $g_\gamma : X_\gamma :=X\times _TP_\gamma  \ra  P_\gamma$ 
with respect to  $p_\gamma $. 
First assume that $f$ admits 
a holomorphic section. 
Then such a section gives rise to 
the tautological line bundle $L_\gamma $ on  $X_\gamma $ 
such that $(L_\gamma)_y, y\in P_\gamma $, is precisely 
the holomorphic line bundle on $X_{\gamma ,y}\cong X_{p_\gamma (y)}$ 
corresponding to $y$ (cf.\ \cite{gr}). 

Then by Corollary \ref{lss} applied to $(g_\gamma,L_\gamma)$,  
for any $k\geq 0$ 
the set  $S_k(L_{\gamma})$ is 
the union of at most countably many analytic subvarieties of $P_\gamma $. 
Then since  $p_\gamma $  is proper, 
by Remmert $\bar{S}_k(L_\gamma ):=p_\gamma(S_k(L_\gamma ))$ also is 
the union of at most countably many analytic subvarieties of $T$, and hence 
so is the union 
\[ A_k:=\bigcup_{\gamma\in \Gamma} \bar{S}_k(L_\gamma )=\{t\in T; \kappa(L_t)\geq k 
\mbox{for some} \ \ L_t\in Pic X_t \} .\]   
On the other hand, 
by (\ref{a}) we have  $a(X_s) \geq  k$  if and only if  $s\in A_k$.  
Namely, $A_k=T_k$ and the proposition is proved in this case.  

When  $f$  does not admit a holomorphic section, we consider the 
pull-back $f_X$  of  $f$  to $X$  via  $f: X \ra  T$, namely  $f_X$ is 
the projection $f_X: X\times _TX \ra  X$  to the second factor.  
Since  $f_X$  admits a canonical holomorphic section, 
we can apply what we have proved above to  $f_X$  and get that 
$X_k:=\{x\in X; a(X_x)\geq k\} $ is 
at most a countable union of analytic subvarieties of  $X$ for all $k$.  
In view of the relation  $f^{-1}(T_k)=X_k$ coming from the definitions of  $T_k$ and $X_k$, 
we deduce the same conclusion for  $T_k$.  Thus the proposition is proved when 
$T$ is simply connected.  

In the general case, 
we consider the pull-back $f_{\tilde{T}}: X\times_T\tilde{T}\ra \tilde{T}$  of  $f$  
to the universal covering $\tilde{T}$ of $T$.  
By what we have proved above, $\tilde{T}_k$ is 
at most a countable union of analytic subvarieties of  $\tilde{T}$ and 
$T_k$ coincides with the image of $\tilde{T}_k$.  
Now if dim $T=1, \tilde{T}_k$ consists of at most countably many points 
and hence so does $T_k$.  Thus the proposition is true also in this case. 
\hfill  
 
\vspace{3 mm} 
{\em Remark 4.1}. 1) If $T$ is not simply connected and dim $T>1$, the above argument 
fails since the image in $T$ of an analytic subvariety of $\tilde{T}$ is not in general 
analytic.  We do not know if the proposition still holds for a general $T$. 

2) The above arguments cannot be applied to the families in Theorem \ref{main} 
since for them the components $P_\gamma$ of Pic $X/T$ are not proper over the base  $T$. 
Indeed, we have $P_\gamma\cong T_\gamma \times \C^*$.

\vspace{3 mm} 

\begin{flushright}
Department of Mathematics\\
Graduate School of Science\\
Osaka University\\
Toyonaka 560-0043 Japan \\ 
{\small E-mail address: fujiki@math.sci.osaka-u.ac.jp} \\ 
Dipartimento di Matematica\\
Universit\'a  Roma Tre\\
00146 Roma, Italy\\
{\small E-mail address: max@mat.uniroma3.it}   
\end{flushright}

\end{document}